\documentclass[12pt]{article}
\begin{document}
\newtheorem{prop}{}[section]
\newtheorem{defi}[prop]{}
\newtheorem{lemma}[prop]{}
\newtheorem{rema}[prop]{}
\def\s{\scriptstyle}
\def\ss{\scriptscriptstyle}
\def\dd{\displaystyle}
\def\a{a}
\def\m{\alpha}
\def\TT{T}
\def\Hinf{ H^{\infty}(\reali^d, \complessi) }
\def\Hn{ H^{n}(\reali^d, \complessi) }
\def\Hm{ H^{m}(\reali^d, \complessi) }
\def\Ha{ H^{\a}(\reali^d, \complessi) }
\def\Ld{L^{2}(\reali^d, \complessi)}
\def\Lpi{L^{p}(\reali^d, \complessi)}
\def\Lq{L^{q}(\reali^d, \complessi)}
\def\Lr{L^{r}(\reali^d, \complessi)}
\def\Knb{K^{best}_n}
\def\sc{\cdot}
\def\k{\mbox{{\tt k}}}
\def\x{\mbox{{\tt x}}}
\def\g{ {\textbf g} }
\def\QQQ{ {\textbf Q} }
\def\AAA{ {\textbf A} }
\def\gr{\mbox{graph}~}
\def\Q{$\mbox{Q}_a$~}
\def\PZ{$\mbox{P}^{0}_a$~}
\def\PZAL{$\mbox{P}^{0}_\alpha$~}
\def\PL{$\mbox{P}^{1/2}_a$~}
\def\PU{$\mbox{P}^{1}_a$~}
\def\PK{$\mbox{P}^{k}_a$~}
\def\PKU{$\mbox{P}^{k+1}_a$~}
\def\PI{$\mbox{P}^{i}_a$~}
\def\Pell{$\mbox{P}^{\ell}_a$~}
\def\PTM{$\mbox{P}^{3/2}_a$~}
\def\AZ{$\mbox{A}^{0}_r$~}
\def\AU{$\mbox{A}^{1}$~}
\def\epsilona{\epsilon^{\scriptscriptstyle{<}}}
\def\epsilonb{\epsilon^{\scriptscriptstyle{>}}}
\def\lgraffa{ \mbox{\Large $\{$ } \hskip -0.2cm}
\def\rgraffa{ \mbox{\Large $\}$ } }
\def\restriction{ \stackrel{\setminus}{~}\!\!\!\!|~}
\def\M{{\scriptscriptstyle{M}}}
\def\Fre{Fr\'echet~}
\def\I{{\mathcal N}}
\def\ap{{\scriptscriptstyle{ap}}}
\def\fiap{\varphi_{\ap}}
\def\BBB{ {\textbf B} }
\def\EEE{ {\textbf E} }
\def\FFF{ {\textbf F} }
\def\GGG{ {\textbf G} }
\def\TTT{ {\textbf T} }
\def\KKK{ {\textbf K} }
\def\FFi{ {\bf \Phi} }
\def\GGam{ {\bf \Gamma} }
\def\sc{ {\scriptstyle{\bullet} }}
\def\ep{\epsilon}
\def\parn{\par\noindent}
\def\teta{M}
\def\elle{L}
\def\ro{\rho}
\def\al{\alpha}
\def\si{\sigma}
\def\be{\beta}
\def\ga{\gamma}
\def\de{\delta}
\def\la{\lambda}
\def\te{\vartheta}
\def\ch{\chi}
\def\et{\eta}
\def\complessi{{\textbf C}}
\def\reali{{\textbf R}}
\def\interi{{\textbf Z}}
\def\naturali{{\textbf N}}
\def\T{{\textbf T}}
\def\T1{{\textbf T}^{1}}
\def\EE{{\mathcal E}}
\def\FF{{\mathcal F}}
\def\GG{{\mathcal G}}
\def\KK{{\mathcal K}}
\def\PP{{\mathcal P}}
\def\QQ{{\mathcal Q}}
\def\J{J}
\def\Np{{\hat{N}}}
\def\Lp{{\hat{L}}}
\def\Jp{{\hat{J}}}
\def\Pp{{\hat{P}}}
\def\Pip{{\hat{\Pi}}}
\def\Vp{{\hat{V}}}
\def\Ep{{\hat{E}}}
\def\Fp{{\hat{F}}}
\def\Gp{{\hat{G}}}
\def\Kp{{\hat{K}}}
\def\Ip{{\hat{I}}}
\def\Tp{{\hat{T}}}
\def\Mp{{\hat{M}}}
\def\Ga{\Gamma}
\def\Si{\Sigma}
\def\Upsi{\Upsilon}
\def\Gag{{\check{\Gamma}}}
\def\Sip{{\hat{\Sigma}}}
\def\Upsig{{\check{\Upsilon}}}
\def\Kg{{\check{K}}}
\def\ellp{{\hat{\ell}}}
\def\j{j}
\def\jp{{\hat{j}}}
\def\BB{{\mathcal B}}
\def\LL{{\mathcal L}}
\def\SS{{\mathcal S}}
\def\DD{{\mathcal D}}
\def\VV{{\mathcal V}}
\def\WW{{\mathcal W}}
\def\OO{{\mathcal O}}
\def\CC{{\mathcal C}}
\def\RR{{\mathcal R}}
\def\AA{{\mathcal A}}
\def\CC{{\mathcal C}}
\def\JJ{{\mathcal J}}
\def\NN{{\mathcal N}}
\def\WW{{\mathcal W}}
\def\HH{{\mathcal H}}
\def\XX{{\mathcal X}}
\def\YY{{\mathcal Y}}
\def\ZZ{{\mathcal Z}}
\def\UU{{\mathcal U}}
\def\CC{{\mathcal C}}
\def\XX{{\mathcal X}}
\def\RR{{\mathcal R}}
\def\cir{{\scriptscriptstyle \circ}}
\def\circa{\thickapprox}
\def\vain{\rightarrow}
\def\parn{\par \noindent}
\def\salto{\vskip 0.2truecm \noindent}
\def\spazio{\vskip 0.5truecm \noindent} 
\def\vs1{\vskip 1cm \noindent}
\def\fine{\hfill $\diamond$ \vskip 0.2cm \noindent}
\newcommand{\rref}[1]{(\ref{#1})}
\def\beq{\begin{equation}}
\def\feq{\end{equation}}
\def\beqq{\begin{eqnarray}}
\def\feqq{\end{eqnarray}}
\def\barray{\begin{array}}
\def\farray{\end{array}}
\makeatletter
\@addtoreset{equation}{section}
\renewcommand{\theequation}{\thesection.\arabic{equation}}
\makeatother
\begin{titlepage}
\begin{center}
{\huge On the constants in some inequalities for 
the Sobolev norms and pointwise product.}
\end{center}
\vspace{1truecm}
\begin{center}
{\large
Carlo Morosi${}^1$, Livio Pizzocchero${}^2$} \\
\vspace{0.5truecm}
${}^1$ Dipartimento di Matematica, Politecnico di
Milano, \\ P.za L. da Vinci 32, I-20133 Milano, Italy \\
e--mail: carmor@mate.polimi.it \\
${}^2$ Dipartimento di Matematica, Universit\`a di Milano\\
Via C. Saldini 50, I-20133 Milano, Italy\\
and Istituto Nazionale di Fisica Nucleare, Sezione di Milano, Italy \\
e--mail: livio.pizzocchero@mat.unimi.it
\end{center}
\vspace{1truecm}
\begin{abstract} \noindent We consider the 
Sobolev norms of the pointwise product of two functions, and estimate 
from above and below the constants appearing in two related inequalities.
\end{abstract}
\vspace{2truecm} \noindent
 \par
\vspace{0.4truecm} \noindent
 \par
\vspace{0.4truecm} \noindent
\textbf{To appear in the Journal of Inequalities and Applications.}
\end{titlepage}
\section{Introduction.} 
\label{intro}
For $d \in \naturali_0 := \naturali \setminus \{ 0 \}$ 
and $n \in [0, +\infty)$, let us
consider the Sobolev (or Bessel potential \cite{Smi} \cite{Maz}) space 
$\Hn$, with the standard norm 
$\|~\|_n$; this is defined setting
\beq   \| f \|_n := \| \sqrt{1 - \Delta}^{~n}~ f \|_{L^2} =
\sqrt{ \int_{\reali^d} d k~ 
(1 + | k |^2)^n ~| (\FF f)(k) |^2 }~, \label{repint} \feq
where $\Delta$ is the Laplacian and 
$\FF$ is the Fourier transform, normalised so that
\beq (\FF f)(k) = {1 \over (2 \pi)^{d/2}} \int_{\reali^d} d x~e^{-i k \sc x}
f(x) \label{ff}~. \feq
(For $n$ integer, $\| f \|_n$ can also be expressed in terms of the 
partial derivatives of $f$ of all orders $\leq n$, see Eq.\rref{inner}).
In the sequel, we consider besides
$n$ another real number $\a$. The following statement is known in the literature: 
\begin{prop}
\label{nm}
\textbf{Proposition.} 
i) Let $n, \a$ be such that $0 \leq n \leq d/2 < \a$. Then,
for each $f \in \Ha$, $g \in \Hn$ it is $f g \in \Hn$; also,
there is a constant $K_{n,\a,d}$ such that 
\beq \| f g \|_n \leq K_{n,\a,d}~ \| f \|_{\a}~ \| g \|_n \label{in1} \feq
for all $f \in \Ha$, $g \in \Hn$. \parn
ii) Let $n > d/2$. Then, for each $f, g \in \Hn$ it is $f g \in \Hn$; 
furthermore, for each $\a$ such that $n \geq \a > d/2$ there 
is a constant $K_{n, \a, d}$ such that 
\beq \| f g \|_n \leq K_{n, \a, d}~ \mbox{Max} 
\Big(\| f \|_{\a}~ \| g \|_n, \| f \|_{n}~ \| g \|_{\a} \Big)
\label{in2} \feq
for all $f, g \in \Hn$.
\fine
\end{prop}
Eq.\rref{in2} makes sense because $H^n \subset H^a$ for $n \geq a$.
For $n$ integer, \rref{in2} follows from the so-called 
"Moser calculus inequality" \cite{Mos} \cite{Maj} \cite{Zei}. 
Both for integer and noninteger $n$ (and $a$), Eq.s \rref{in1} and
\rref{in2} can be obtained 
specializing to the classical Sobolev spaces a more 
general result on the Triebel-Lizorkin spaces: see, e.g., \cite{Run}. 
The inequalities in Prop.\ref{nm}
have interesting applications to nonlinear PDE's; in particular, they can be 
employed to derive "tame" estimates (in the sense of the Nash-Moser 
theory) in the case of polynomial nonlinearities 
\cite{Ham} \cite{MP} (in these applications, one generally 
asks $n$ to be an arbitrary integer, and $\a$ the 
smallest integer $> d/2$, i.e., $\a = [d/2] + 1$).
\parn
From here to the end of the paper, \textsl{we intend that} $K_{n, \a, d}$ 
\textsl{is the sharp (i.e., the minimum) constant} satisfying for 
all $f, g$ the inequality \rref{in1} if $0 \leq n \leq d/2 < \a$, 
or the inequality \rref{in2} if $n \geq \a > d/2$. \parn
Neither in the quoted references, nor in any other 
of our knowledge, these constants are estimated; this is a nontrivial task, 
because essentially one must evaluate the quantity
$\| \sqrt{1 - \Delta}^{~n}~ (f g) \|^2_{L^2} =$
$\int_{\reali^d} d k~ 
(1 + | k |^2)^n ~| \FF (f g)(k) |^2$, quartic in 
the pair $(f, g)$. A direct variational approach seems to be 
very difficult, even for $n$ and $a$ integers (to say the least:
the formal stationarity condition of the functional to be maximised
for finding $K_{n,a,d}$  
yields a system of 
cubic PDE's of order $2~ \mbox{Max}(n, a)$ for the pair $(f, g)$; 
one can hope to treat it for special values of $n$,$a$ at most). \parn
On the other hand explicit estimates on $K_{n,a,d}$, even 
non optimal, holding for all values of 
$n, a$ (and $d$) would be useful for the strictly quantitative aspects of the
previously mentioned applications; the aim of the present article 
is just to give some results of this sort.
The paper is organized as follows. \parn
In the next Section we describe 
the results to be derived in the subsequent ones; these are, essentially,
the upper bounds on $K_{n,a,d}$ coming from a general argument,
and the lower bounds obtained by 
substituting convenient classes of 
trial functions in Eq.s \rref{in1} \rref{in2}. These bounds allow to 
estimate $K_{n,a,d}$ from above and below, for general values of 
$n,a,d$. The explicit numerical values of the bounds are reported for a number 
of cases with low $n,a$ and $d$; for large $n$, our
estimates have the 
form $\mbox{const.} 2^n/(n+a)^{a/2 + d/4} \stackrel{<}{\sim} K_{n,a,d} 
\stackrel{<}{\sim} \mbox{const.} 2^n$. \parn
In Sect.\ref{sobolev} our 
normalizations and notational conventions on Sobolev spaces and 
their norms are fixed; also, we report an estimate \cite{art} on the 
constants in the classical imbedding inequalities of Sobolev into $L^r$ spaces, 
to be employed later on. 
\parn
In Sect.\ref{oper} we define
some auxiliary nonlinear operators 
$f \mapsto \DD_{n}(f)$ on Sobolev spaces, 
that will be a basic tool to prove our upper bounds on 
the constants $K_{n,a,d}$; as a first step towards this goal, 
we infer a non conventional "Leibnitz" inequality for $\DD_n(f g)$, 
where $f$, $g$ are two functions. \parn
In Sect.\ref{deriva} we derive some estimates on pointwise 
products of the form $\DD_{l}(f) \DD_{m}(g)$; in the 
subsequent Sect.\ref{above}, these estimates are employed to obtain a
new proof of Prop.\ref{nm}, and to
to derive our upper bounds on $K_{n, a, d}$. \parn
In the final Sect.s \ref{ground}, \ref{bessel} and  \ref{below} 
we derive the already 
mentioned, different kinds of lower bounds on $K_{n,a,d}$, corresponding to 
different trial functions. 
\section{Description of the results.}
\label{main} 
Many of our bounds will be given in terms of 
the usual $\Gamma$ function, and of the function $E$ defined by
\beq E(s) := s^s \quad \mbox{for $s \in (0, + \infty)$}~, \qquad
E(0) := 1~. \label{defe} \feq
We will employ the coefficients
\beq S_{\a, d} := {1 \over (4 \pi)^{d/4}} \sqrt{ 
{\Gamma(\a - d/2) \over \Gamma(\a)} }~, \label{sad} \feq
\beq E_{\ell, \a, d} := 
\left( { E\left({\displaystyle{ {\ell \over 2 \a} }}\right)
E\left({\displaystyle{ {1 \over 2}  - {\ell \over 2 \a} }}\right) \over 
E\left({\displaystyle{ {1 \over 2}  + {\ell \over 2 \a} }}\right) 
E\left({\displaystyle{ 1 - {\ell \over 2 \a} }}\right) }  \right)^{d/2} 
\label{elad} \feq
($ \a > d/2$, $~ 0 \leq \ell \leq \a$); 
as explained in the sequel, these can be interpreted 
in terms of the imbedding inequalities of Sobolev into $L^{r}$ spaces. 
One finds by elementary means 
that the function $\ell \in [0, \a] \mapsto E_{\ell, \a, d}$
attains its maximum at $\ell = 0,~ \a$ and its minimum at $\ell= \a/2$; 
these are, respectively,
\beq E_{0,\a, d} = E_{\a, \a, d} = 1~, 
\qquad E_{\a/2, \a, d} = \left(16/27\right)^{d/4}~. \feq
Our estimates will rely on some combinations of the 
above constants with the binomial coefficients; to deal with Sobolev 
spaces of noninteger order, certain generalised binomial coefficients
will be necessary. For each $n \in [0, +\infty)$,  
we will put
\beq n_{+} := \mbox{Min} \left\{ m \in \naturali~|~m \geq n \right\}~; \feq
(clearly, $n_{+}$ is an integer approximation of $n$ from above, not 
to be confused with the standard integer part $[n] :=$ 
maximum integer $\leq n$). \parn
We will employ the "lattice"
(with initial point $0$ and final point $n$)
\beq \Lambda(n) := \{ j~{n \over n_{+}}~|~j=0,..., n_{+} \} 
\label{lat} \feq
and the coefficients 
\beq \left( \barray{cc} n \\ \ell \farray \right)_{+} := 
\left( \barray{cc} n_{+}  \\ j \farray \right) \qquad
\mbox{for $\ell = j~ \displaystyle{ {n \over n_{+}} } 
\in \Lambda(n)$}~. \label{bin} \feq
In the above formula, $(~)$ denotes the usual binomial 
coefficients; also, it is understood that $n_{+}/n := 1$ if $n=0$.
For any $n \in \naturali$, it is $n_{+}= n$, $\Lambda(n) =  
\{ 0, 1, ..., n \}$, and 
$\displaystyle{\left( \barray{cc} n \\ \ell \farray \right)_{+}}$
are the usual binomial coefficients 
$\displaystyle{\left( \barray{cc} n \\ \ell \farray \right)}$. In general,
we have
\beq \sum_{\ell \in \Lambda(n)} 
\left( \barray{cc} n \\ \ell \farray \right)_{+} = 2^{n_{+}}~. \label{duen} 
\feq
\vskip 0.2cm\noindent
\textbf{General upper bounds on} $\mbox{\boldmath $K_{n, \a, d}$}$~.
These will be obtained by combining the Sobolev 
imbedding inequalities of $H^n$ into $L^r$ spaces (Sect.\ref{sobolev})
with a "Leibnitz" inequality for the $H^n$ norms of products 
(Sect.\ref{oper}); the final result will be the following.
\begin{prop}
\label{upper}
\textbf{Proposition.} Let $d \in \naturali_0$ and either 
$0 \leq n \leq d/2 < \a$ or $n \geq \a > d/2$; then the sharp 
constant in the inequality \rref{in1} or \rref{in2} admits the 
upper bound
\beq K_{n, \a, d} \leq S_{\a, d}~\sum_{\ell \in \Lambda(n)}~
\left( \barray{cc} n \\ \ell \farray \right)_{+}~E_{n, \ell, \a, d}~, 
\label{ub} \feq
where the coefficients in the r.h.s. are defined as follows. \parn
In the case $0 \leq n \leq d/2 < \a$, we put
\beq E_{n, \ell, \a, d} := E_{\ell, \a, d}~, \label{knelad1} \feq
while for $n \geq \a > d/2$ we put
\beq E_{n, \ell, \a, d} := \left\{ \barray{lll} E_{\ell, \a, d} && 
\mbox{if~ $0 \leq \ell < \a/2$}~, \\ 
(16/27)^{d/4}  && \mbox{if~ $\a/2 \leq \ell \leq n - \a/2$}~, \\
E_{n- \ell, \a, d} && \mbox{if~ $n - \a/2 < \ell \leq n$~. } 
\farray \right. \label{knelad2} \feq
Eq.s (\ref{ub}-\ref{knelad2}) imply weaker bounds: for
all $n, \a$, they give
\beq K_{n, \a, d} \leq S_{\a, d}~2^{n_+} \label{wub1} \feq
and, for $n \geq a > d/2$, they imply 
\beq K_{n, \a, d} \leq \left(16/27 \right)^{d/4} S_{\a, d}~
u_{n, \a, d}~2^{n_{+}}~, 
\label{wub2} \feq
\beq u_{n, \a, d} := 1 + { (27/16)^{d/4} - 1 \over 2^{n_+ - a_n} } 
\left( \barray{cc} n_+ \\ n_{+} - a_n + 1 \farray \right)~, 
\qquad a_{n} := \left( {n_+ \over n}~{a \over 2} \right)_{+}~. 
\label{defna} \feq
(Note that, by construction, 
$\lim_{n \vain +\infty} u_{n, \a, d} = 1$
for fixed $\a$ and $d$). 
\fine
\end{prop}
\textbf{General method to derive lower bounds on} 
$\mbox{\boldmath $K_{n, \a, d}$}$~. Of course, Eq.s \rref{in1} \rref{in2}
imply
\beq K_{n, a, d} \geq { \| f g \|_{n} \over \| f \|_a \| g \|_n } 
\label{ofc1} \feq
for $0 \leq n \leq d/2 < a$ and any nonzero 
$f \in \Ha$, $g \in \Hn$~, and
\beq K_{n, a, d} \geq { \| f g \|_{n} \over \mbox{Max} 
(\| f \|_a \| g \|_n, \| f \|_n \| g \|_{a}) } \label{ofc2} \feq
for $n \geq a > d/2$ and any nonzero $f, g \in \Hn$~. \parn
All the results presented hereafter will be derived inserting 
convenient trial functions $f$, $g$ into Eq.s \rref{ofc1} \rref{ofc2}.
In particular, we will use the functions (or some rescaled variants)
\beq f_{n, d} := \FF^{-1}\left(1 \over (1 + | \k |^2)^n \right) \qquad 
(n > 0)~, \label{fnd} \feq
where $\FF$ is again the Fourier transform, and
$| \k |$ denotes the map 
$(k_1,...,k_d) \in \reali^d $ $\mapsto \sqrt{{k_1}^2 + ... + {k_d}^2}$. 
From the above definition, it is clear that $f_{n, d} \in \Hn$ if $n > d/2$.
By a known formula for radially simmetric
Fourier transforms \cite{Boc}, we have \cite{Smi} \cite{Maz}
\beq f_{n, d} = {| \x |^{n - d/2} \over 2^{n - 1} \Gamma(n)}~
K_{n - d/2}(| \x |)~, \label{gemac} \feq
where $| \x |$ is the function $x \mapsto \sqrt{{x_1}^2 + ... {x_d}^2}$ 
and $K_{(~)}$ are the
modified Bessel functions of the third kind, or Macdonald functions, see, 
e.g., \cite{Wat}.  
\parn
Another family of trial functions, useful for our purposes, is made 
of the functions
\beq x = (x_1, ..., x_d) \mapsto f_{p, \sigma, d}(x) := 
e^{i p x_1}~e^{- (\sigma/ 2) | x |^2} \qquad (~p, \sigma 
\in (0, +\infty)~) \label{fg} \feq
(with $| x | := \sqrt{{x_1}^2 + ... + {x_d}^2}$); here, the Fourier 
character $e^{i p x_1}$ is regularised at infinity by the rapidly 
vanishing, Gaussian factor $e^{- (\sigma/ 2) | x |^2}$. 
We shall mainly deal 
with these functions in the limit of small $\sigma$ and large $p$. \parn
Let us present three kinds of lower bounds obtained from the above
trial functions; 
for fixed $d$ and $a$, these are 
interesting for $n$ very low, $n$ close to $a$ and $n$ large, 
respectively. 
\vskip 0.2cm \noindent
\textbf{"Ground level" lower bounds on} $\mbox{\boldmath $K_{n, \a, d}$}$~. 
First of all one can show that, in any case, 
$K_{n, a, d}$ is bounded from below by a constant independent of $n$.
Either in Eq.\rref{ofc1} or in
\rref{ofc2}, we insert the trial functions $f := f_{a, d}$ (see Eq.s 
\rref{fnd} \rref{gemac}) and 
$g := $ a smooth approximant of the Dirac $\delta$ distribution. 
This yields the following estimate.
\begin{prop}
\label{gen}
\textbf{Proposition.} Let $d \in \naturali_0$ 
and either $0 \leq n \leq d/2 < \a$ or $n \geq \a > d/2$; then
\beq K_{n, a, d} \geq S_{a, d}~. \label{bsad} \feq 
\fine
\end{prop}
The constant $S_{a, d}$ is the same appearing in the upper bounds
\rref{ub}. In particular, for $n=0$ the
upper and lower bounds \rref{ub} \rref{bsad} coincide, allowing 
to individuate the sharp constant; in other terms, we have
\begin{prop}
\label{cogen}
\textbf{Corollary.} For $d \in \naturali_0$, $\a > d/2$ and 
$n=0$ the sharp constant in the inequality \rref{in1} is
\beq K_{0, a, d} = S_{a, d}~. \label{usad} \feq 
\fine
\end{prop}
For $n$ low, the upper and lower bounds \rref{ub} \rref{bsad} are not far,
thus confining the sharp constant to a fairly small interval; for 
example, if $n=1 \leq d/2$ this interval is described by the inequalities
\beq S_{a, d} \leq K_{1, a, d} \leq S_{a, d} ~(1 + E_{1,a,d})~. 
\label{cogen1} \feq 
\vskip 0.2cm \noindent
\textbf{"Bessel" lower bounds on } 
$\mbox{\boldmath $K_{n, \a, d}$}$~. We assume $n \geq a > d/2$, and 
insert into Eq.\rref{ofc2} the functions 
$f := g := f_{\la, n, d}$,
where
\beq f_{\la, n, d} := 
{1 \over \la^d} \FF^{-1}\left(1 \over (1 + | \k |^2/\la^2)^n \right)
\qquad (\la > 0)~; \label{lagemac} \feq
these come from rescaling by a factor $\la$ of the function $f_{n, d}$ in 
Eq.\rref{fnd}, i.e., $f_{\la, n, d}(x) = f_{n, d}(\la x)$. 
We write down Eq.\rref{ofc2} for these functions and maximise
w.r.t. $\lambda$; the conclusions stemming from this analysis
can be summarised as follows.
\begin{prop}
\label{lalem}
\textbf{Proposition.} For $d \in \naturali_0$, $n \geq a > d/2$ the sharp
constant in Eq.\rref{in2} is such that
\beq K_{n, a, d} \geq \mbox{Sup}_{\lambda > 0} 
{ \|  f_{\la, n, d}^2 \|_n \over 
\| f_{\la, n, d} \|_{a} ~\| f_{\la, n, d} \|_n }~.
\label{dare}
\feq
The norms in the r.h.s. of Eq.\rref{dare} can be expressed in terms 
of hypergeometric functions (or one-dimensional integrals of them).
\fine
\end{prop}
The explicit expressions of the above norms 
will be given in Sect.\ref{bessel}. Experimentally, the 
lower bounds \rref{dare} are not too far from the upper bounds 
\rref{ub} when $n$ and $a$ are fairly close to $d/2$. \parn
\vskip 0.2cm \noindent
\textbf{Numerical examples corresponding to the upper 
bounds \rref{ub} and to the "ground level" and "Bessel" lower bounds.} 
The numerical values reported hereafter have been
obtained from the (analytical) estimates
mentioned previously by means of the MATHEMATICA package. 
We consider the cases $d=1,2,3$, $a = [d/2] + 1$ 
and $n$ integer, $0 \leq n \leq a$. \parn
Let $d=1$, $a=1$; for $n=0$ we apply Eq.\rref{usad}, whereas 
for $n=1$ we use Eq.s \rref{dare}\rref{ub}; this yields the 
estimates
\beq 0.71 < K_{0,1,1} = 1/\sqrt{2} < 0.72~, \feq
\beq 0.84 < K_{1,1,1} < 1.42~. \label{just} \feq
Let $d=2$, $a=2$; using Eq.\rref{usad} for $n=0$, 
Eq.\rref{cogen1} for $n=1$ and Eq.s \rref{dare} \rref{ub} for $n=2$ 
we get, respectively, 
\beq 0.27 < K_{0,2,2} = 1/(2 \sqrt{\pi}) < 0.28~, \feq
\beq 0.27 < K_{1,2,2} < 0.50~, \feq
\beq 0.36 < K_{2,2,2} < 1.00~. \label{ju36} \feq
Finally, let $d=3$, $a=2$; using again Eq.\rref{usad} for $n=0$, 
Eq.\rref{cogen1} for $n=1$ and Eq.s \rref{dare} \rref{ub} for $n=2$ 
we get, respectively,
\beq 0.19 < K_{0,2,3} = 1 / (2 \sqrt{2 \pi}) < 0.20~, \feq
\beq 0.19 < K_{1,2,3} < 0.34~, \feq
\beq 0.24 < K_{2,2,3} < 0.67~. \label{ju24} \feq
\vskip 0.2cm \noindent
\textbf{"Fourier" lower bounds on} 
$\mbox{\boldmath $K_{n, \a, d}$}$~. We 
insert into Eq.\rref{ofc1} or \rref{ofc2} the trial functions
$f := g := f_{p, \sigma, d}$,~where $f_{p, \sigma, d}(x) := 
e^{i p x_1}~e^{- (\sigma/ 2) | x |^2}$ as in Eq.
\rref{fg}. As anticipated, 
the Gaussian factor $e^{-(\sigma/2) |x|^2}$ is used only to 
regularise at infinity 
the Fourier character $e^{i p x_1}$; it will be ultimately
taken as close as possible to unity, setting $\sigma$ close to zero. 
The norms of the functions \rref{fg} are 
evaluated in Sect.\ref{below}; in the same Section, we will 
subsequently choose $p = \mbox{const.} \sqrt{n +a}$, 
$\sigma = \mbox{const.}/\sqrt{n+a}$. The estimate arising in this way 
is interesting for $n$ great, but can be nominally written for small $n$ 
too; it can be expressed as follows. 
\begin{prop}
\label{lower}
\textbf{Proposition.} Let $d \in \naturali_0$, and either 
$1/2 \leq n \leq d/2 < \a$ or $n \geq \a > d/2$; then, the sharp 
constant in Eq.\rref{in1} or \rref{in2} admits the lower 
bounds
\beq K_{n, \a, d} \geq R_{\a, d}~v_{n, \a, d}~
{2^n \over (n + \a)^{\a/2 + d/4}}~, \label{lb} \feq
\beq R_{\a, d} := {e^{-\a/2} \over (2 \pi)^{d/4} }~
\sqrt{ E(d/2) E(\a - d/2) }~, \label{rad} \feq
\beq v_{n, \a, d} :=  
\left( 1 - {d \over 2 (n + a)} + 
{d^2 a n \over 4(n + \a)^4} \right)^{d/4}~
e^{\dd{ {(4 \a - d) d n + 2 d \a^2 
\over 8 (n + \a)^2 - 4 d n }  - 
{d \a^2 \over 4 (n + \a)^2 - 2 d \a} } }~. \label{vnad} \feq
(Note that $\lim_{n \vain +\infty} v_{n, \a, d} = 1 $
for fixed $\a$ and $d$). \parn 
For $n \geq \a > d/2$, Eq.\rref{lb} implies the weaker bound
\beq K_{n, \a, d} \geq 
R_{\a, d}~\left(1 - {d \over 2 \a} \right)^{d/4}~
{2^n \over (n + \a)^{\a/2 + d/4}}~. \label{wlb} \feq
\fine
\end{prop}
\vskip 0.2cm \noindent
\textbf{The} $\mbox{\boldmath $n \vain + \infty$}$ 
\textbf{limit for the constants} $\mbox{\boldmath $K_{n,a,d}$}$. 
From Prop.s \ref{upper} and \ref{lower} we see that $K_{n, a, d}$ 
has upper and lower bounds behaving essentially as $2^n$ and
$2^n/n^{a/2 + d/4}$, respectively. 
These yield upper and lower bounds for the ratio 
$(1/n) \log_{2} K_{n, \a, d}$, both converging to  
$1$ for $n \vain \infty$. So, Prop.s \ref{upper} and 
\ref{lower} imply
\begin{prop}
\label{corolla}
\textbf{Corollary.} For fixed $\a$ and $d$, it is 
\beq {\rm{lim}}_{n \vain +\infty}~
{\log_{2} K_{n, \a, d} \over n} = 1~. \feq
\fine
\end{prop}
\noindent
\section{Some facts on Sobolev spaces.}
\label{sobolev}
\textbf{Notations.} We stick to the previous paper \cite{art}. 
As usually: $d \in \naturali_0 = \naturali \setminus \{ 0 \}$ is 
a fixed space dimension; the running variable in $\reali^d$ is 
$x = (x_1, ..., x_d)$, and $k = (k_1, ..., k_d)$ when using the 
Fourier transform; we write $| \x |$ 
for the function $(x_1, ..., x_d) \mapsto \sqrt{{x_1}^2 + ... + {x_d}^2}$, and
intend $| \k |$ similarly.
We denote with 
\beq \FF, \FF^{-1} : \SS'(\reali^d, \complessi) \vain \SS'(\reali^d, 
\complessi) \feq
the Fourier transform of tempered distributions 
and its inverse, choosing normalizations so that \rref{ff} holds
(whenever the integral makes sense, say, for
$f$ in $L^1(\reali^d, \complessi)$~). 
The restriction of $\FF$ to $\Ld$, with the standard
inner product and the associated norm $\|~\|_{L^2}$, is 
a Hilbertian isomorphism. \parn
Consider a real number $n \geq 0$, and introduce the operator 
\beq \SS'(\reali^d, \complessi) \vain \SS'(\reali^d, \complessi)~, 
\qquad f \mapsto \sqrt{1 - \Delta}^{~n}~ f 
:= \FF^{-1} \left( \sqrt{1 + | \k |^2}^{~n} \FF f \right) \label{lap} \feq
(if $n$ is an even integer, this is a power of 
$1$ minus the distributional Laplacian $\Delta$, in the elementary
sense). The $n$-th order Sobolev
(or Bessel potential) space of $L^2$ type, 
denoted with $\Hn$, and its norm $\|~\|_n$ are defined as follows:
\beq \Hn := \lgraffa f \in \SS'(\reali^d, \complessi)~\Big\vert~
\sqrt{1 - \Delta}^{~n} f \in \Ld~ \rgraffa= \label{incid} \feq
$$ = \lgraffa f \in \SS'(\reali^d, \complessi)~\Big \vert~
\sqrt{1 + | \k |^2}^{~n}  \FF f \in \Ld \rgraffa~, $$
\beq \| f \|_{n} := \| \sqrt{1 - \Delta}^{~n}~ f \|_{L^2} = 
\|~ \sqrt{1 + | \k |^2}^{~n}~\FF f~ \|_{L^2}~. \label{repfur} \feq
Of course, if $n \leq n'$, it is 
$H^{n'}(\reali^d, \complessi) \subset \Hn$ and 
$\|~\|_{n} \leq \|~\|_{n'}$; 
also, $H^{0} = L^2$. \parn
The connections between $H^n$ spaces and Bessel functions are known 
after \cite{Smi}. For completeness, let us write down the expression
of $\|~ \|_n$ in terms of distributional derivatives, when $n$ is 
integer; in this case, we have 
\beq \| f \|_n = \sqrt{ \sum_{m=0}^n \left( \barray{cc} n \\ m \farray 
\right) \sum_{\alpha \in \naturali^d, | \alpha | = m}~{m! \over \alpha!}~~ 
\| \partial^{\alpha} f \|^{2}_{L^2}~}~, \label{inner} \feq
where, for each $\alpha = (\alpha_1, ..., \alpha_d)$, it 
is intended that $\partial^{\alpha} := {\partial_1}^{\alpha_1} ... 
{\partial_d}^{\alpha_d}$,~ $| \alpha | := \alpha_1 + ... + \alpha_d$,~
$\alpha! := \alpha_1! ... \alpha_d!$~. Eq.\rref{inner} can be derived
with appropriate manipulations from Eq. \rref{repfur}; the 
sums over $m$ and $\alpha$ come essentially from the expansion
of $(1 - \Delta)^n$ or $(1 + | \k |^2)^n$. 
\vskip 0.2cm \noindent
\textbf{Imbedding and interpolation 
inequalities.} For real $r \geq 2$ or $r = \infty$, we consider
the space $L^r(\reali^d, \complessi)$  and its 
norm $\|~\|_{L^r}$. For each real $n \geq 0$, 
$\|~\|_n$ is again the Sobolev norm \rref{repfur}; 
$E(~)$ is the function of Eq.\rref{defe}.
\parn
The imbedding inequalities of $H^n$ into $L^{r}$ spaces are well 
known; in this paper, we will use the following estimate on the imbedding 
constants \cite{art}.
\begin{prop}
\label{sob}
\textbf{Proposition} (Imbedding inequality). 
Let $d \in \naturali_0$, and either 
\beq n=0,~ r=2~, \qquad \mbox{or} \qquad 0 < n < {d / 2},~ 2 \leq r 
< {d \over d/2 - n}~, \qquad \mbox{or} \feq
$$ n = {d / 2},~2 \leq r < \infty~, \qquad \mbox{or} 
\qquad n > {d / 2},~2 \leq r \leq \infty~. $$
Then $\Hn \subset \Lr$, and for each $f \in \Hn$ it is 
\beq \| f \|_{L^r} \leq S_{r, n, d}~ \| f \|_{n}~, \label{sobo} \feq
\beq S_{r, n, d} :=
{1 \over (4 \pi)^{d /4 - d/(2 r)} }~
\left( { \Gamma \left( {\displaystyle{ {n \over 1 - 2/r} - {d \over 2} }} \right) 
\over \Gamma \left( {\displaystyle{ {n \over 1 - 2/r} }} \right) } 
\right)^{1/2 - 1/r}
\left({E(1/r) \over E(1 - 1/r)} \right)^{d/2}~~(r \neq 2, \infty),
\label{srd} \feq
\beq S_{2, n, d} := 1~, \qquad S_{\infty,n,d} := {1 \over (4 \pi)^{d/4}}~
\left( \Gamma(n-d/2) \over \Gamma(n) \right)^{1/2}~.
\label{srdi} \feq
$S_{r,n,d}$ is the sharp (i.e., minimum) imbedding constant 
for $n \geq 0$, $r=2$, and $n > d/2$, $r=\infty$; concerning 
the latter case, the equality $\| f \|_{L^{\infty}} = 
S_{\infty, n, d} \| f \|_n$ holds if $f = f_{n, d} =$ the function 
in Eq.s \rref{fnd} \rref{gemac}. 
\end{prop}
\textbf{Remarks.} i) The constant denoted with $S_{a, d}$ in Eq.\rref{sad}
is just $S_{\infty, a,d}$ with the notations of the above Prop.\ref{sob}. \parn
ii) The estimates of Prop.\ref{sob} on the imbedding constants come 
from application of the Hausdorff-Young and H\"older inequalities to the 
Fourier transform, as shown in \cite{art}; they were previously 
written, for particular cases, in a number of works mentioned
therein. \parn
iii)  For $0 < n < d/2$, 
the imbedding inequality of $\Hn$ into 
$L^r(\reali^d, \complessi)$ still holds in the limit case 
$r=d/(d/2 - n)$, not covered by the previous Proposition. \parn
An analysis of the reliability of the estimates \rref{srd} \rref{srdi} 
was made in \cite{art}; here, the statement that 
$S_{r,n,d}$ is the sharp 
imbedding constant for $r=2$ and 
$n > d/2$, $r=\infty$ was completed showing that
$S_{r,n,d}$ is generally very close to the (unknown) sharp 
constant for $n > d/2$ and arbitrary 
$r$ in $(2, \infty)$. \parn
For 
$0 < n \leq d/2$, one could derive different estimates on the imbedding 
constants using the Hardy-Littlewood-Sobolev inequality, with the 
method indicated in \cite{Miz} and the expression determined in \cite{Lie1} for 
the sharp Hardy-Littlewood-Sobolev constant. (The imbedding constants 
derived in this way become, in the particular cases 
$r=d/(d/2 - n)$, $n=1$ and $2$, the constants determined 
in papers \cite{Aub} \cite{Tal} prior to \cite{Lie1}, and 
in \cite{Wan}, respectively; these are sharp for the inequalities 
considered therein,
strictly related to the above indicated cases of \rref{sobo}). \parn
Nevertheless, all the numerical experiments we performed using 
the alternative, Hardy-Littlewood-Sobolev estimates on the 
imbedding constants yielded no 
essential improvement for the main purpose of the present paper, i.e.,
estimating the constants $K_{n,a,d}$ in the product inequalities \rref{in1} 
\rref{in2}. For this reason, whenever the imbedding constants
will be needed in the sequel, we will always stick to the 
result of Prop.\ref{sob}.  
\fine
The second 
inequality we need in the sequel is the following, known in the literature.
\begin{prop}
\label{segmento}
\textbf{Proposition} (Multiplicative interpolation inequality).
Let $(b, b'), (i, i'), (c, c') \in [0, \infty)^{2}$; assume 
$(i, i')$ to lie on the segment of extremes $(b, b')$ and $(c, c')$, i.e.,
$(i, i') = (1 - t) (b, b') + t (c, c')$ with $t \in [0, 1]$. \parn
Then, for each $f \in H^{Max(b, c)}(\reali^d, \complessi)$, 
$g \in H^{Max(b', c')}(\reali^d, \complessi)$ it is
\beq \| f \|_{i}~ \| g \|_{i'} 
\leq \mbox{Max}\Big(\| f \|_{b}~ \| g \|_{b'}, 
\| f \|_{c}~ \| g \|_{c'} \Big)~. \label{int2} 
\feq
\end{prop}
\textbf{Proof.} We have the inequalities 
(\textsl{\`a la} Browder-Ehrling-Gagliardo-Niremberg \cite{Ada})
\beq \| f \|_{i} \leq \| f \|_{b}^{1 - t}~\| f \|_{c}^{t}~, 
\qquad \| g \|_{i'} \leq \| g \|_{b'}^{1 - t}~\| g \|_{c'}^{t}
\label{int1} \feq
(understanding that $0^t := 0$; these can be derived from the Fourier 
representation \rref{repfur} of the Sobolev norms and H\"older's inequality). 
The thesis \rref{int2} follows multiplying, and 
recalling the elementary inequality 
$v^{1 - t} w^{t} \leq \mbox{Max}(v, w)$ for  $v,w \in [0, + \infty)$ 
and $t \in [0, 1]$. (Similar arguments were employed in \cite{Mos}, and 
in \cite{Ham} for the Sup-norms). \fine
\section{Operators $\mbox{\boldmath $\DD_n$}$. "Leibnitz" inequality.}
\label{oper}
As anticipated, these maps will be basic in our strategy to estimate 
the pointwise product of two functions. \parn
\begin{defi}
\label{defidn}
\textbf{Definition.} For each $n \in [0, +\infty)$, we put
\beq \DD_n : \Hn \vain \Ld~, \qquad f \mapsto \DD_n(f) := 
\FF^{-1} \left( {\sqrt{1 + | \k |^2}}^{~n}  | \FF f | \right)~. 
\label{ddn} \feq 
\fine
\end{defi}
The notation $| \FF f |$ stands for the function 
$\reali^d \vain [0, +\infty)$,~ $k \mapsto | (\FF f)(k) |$;
due to the presence of the modulus $| ~|$, the map $\DD_n$ 
is homogeneous but nonlinear, and differs (on 
its domain) from the linear operator
$\sqrt{1 - \Delta}^{~n}$ of Eq.\rref{lap}. However, comparing with 
the Fourier representations of Sobolev spaces and their norms in Eq.s 
\rref{incid} \rref{repfur}, we see that 
the definition of $\DD_n$ is well posed, and that
\beq \| f \|_n = \| \DD_n(f) \|_{L^2} \qquad \mbox{for $n \in [0, +\infty)$,
$f \in \Hn$}~; \feq
\beq 
\DD_{\ell} \left(H^{\ell + m}(\reali^d, \complessi) \right) 
\subset H^{m}(\reali^d, \complessi), \qquad 
\| \DD_{\ell}(f) \|_{m} = \| f \|_{\ell + m}~, \feq
$$ \DD_{m} (\DD_\ell(f)) = \DD_{\ell + m}(f) 
\qquad \mbox{for $\ell, m \in [0, +\infty),~ 
f \in H^{\ell + m}(\reali^d, \complessi)$}~. $$
Let us consider the action of $\DD_n$
on the pointwise product of two functions $f$ and $g$, and
obtain an inequality concerning $\| \DD_n(f g) \|_{L^2} = 
\| f g \|_n$; this is stated in the following Proposition, to be proved 
after two Lemmas. 
\begin{prop}
\label{disug}
\textbf{Proposition} ("Leibnitz" inequality). 
Let $n \in [0, +\infty)$, and $f, g \in \Hn$ be such that
$\DD_{\ell}(f) \DD_{n- \ell}(g)$ $\in \Ld$ for 
each $\ell \in \Lambda(n)$. Then $f g \in \Hn$, and 
\beq \| f g \|_n  \leq \sum_{\ell \in \Lambda(n)} 
\left( \barray{cc} n \\ \ell \farray 
\right)_{+} \| \DD_{\ell}(f) \DD_{n - \ell}(g) \|_{L^2}~. 
\label{apart} \feq
\fine
\end{prop}
\textbf{Remarks.} i) Recall that $\Lambda(n)$, 
$\left( \barray{cc} n \\ \ell \farray 
\right)_{+}$ are defined by Eq.s \rref{lat} \rref{bin}. \parn
ii) For $n$ integer, one could express $\| f g \|_n$ via Eq.\rref{inner}, 
and the partial derivatives $\partial^{\alpha} (f g)$ appearing therein
in terms of products $\partial^{\lambda} f~ \partial^{\mu} g$ 
($\lambda, \mu$ multiindices),
with the usual Leibnitz rule. Apart from working also for noninteger $n$, 
the estimate \rref{apart} is more efficient for 
evaluating the constants $K_{n,a,d}$ in Eq.s \rref{in1} \rref{in2}.
\fine
Here are the two Lemmas to be employed for proving Prop.\ref{disug}.
\begin{lemma} 
\label{newton}
\textbf{Lemma.} For all real $a, b > 0$, $n \geq 0$ it is
\beq (a + b)^n \leq \sum_{\ell \in \Lambda(n)}~
\left( \barray{cc} n \\ \ell \farray 
\right)_{+} a^{\ell}~b^{n - \ell}~. \label{neweq} \feq
\end{lemma}
\textbf{Proof.} It is obvious if $n=0$. Let $n > 0$; 
being $0 < n/n_{+} \leq 1$, we have
$$ (a + b)^n = \left( (a + b)^{n/n_+} \right)^{n_+} 
\leq \left(a^{n/n_+} + b^{n/n_+}\right)^{n_+}~;  $$
the thesis follows expanding the last expression with
the standard binomial formula for integer exponent $n_+$. 
\fine
\begin{lemma} 
\label{vettori}
\textbf{Lemma.} For all $\xi, \eta \in \reali^d$, it is
\beq \sqrt{1 + | \xi + \eta |^2} < \sqrt{1 + | \xi |^2} + 
\sqrt{1 + | \eta |^2}~. \label{veteq} \feq
\end{lemma}
\textbf{Proof.} It is: $1 + | \xi + \eta |^2 \leq 
1 + | \xi |^2 + | \eta |^2 + 2~ | \xi |~| \eta | < $
$(1 + | \xi |^2) + (1 + | \eta |^2) + 2 \sqrt{1 + | \xi |^2}~
\sqrt{1 + | \eta |^2}$. \fine
The last tool we need to prove Prop.\ref{disug} is the 
convolution product. 
Let us write $F \ast G$ for the convolution of two distributions 
$F, G \in \DD'(\reali^d, \complessi)$, whenever this exists as an 
element of $\DD'(\reali^d, \complessi)$ (see, e.g., \cite{Vla}). We have
\beq (F \ast G)(k) = \int_{\reali^d} 
d h~F(k - h) G(h)  \feq  
if $F, G \in L^1_{loc}(\reali^d, \complessi)$ and the integral 
in the r.h.s. exists, defining an $L^1_{loc}$ function of $k$; 
these conditions are satisfied, in particular, if $F, G \in \Ld$, 
which is the case considered in the forthcoming proof.
With the chosen normalizations for $\FF$ and $\ast$, we have
$\FF (f g) = (2 \pi)^{-d/2}~(\FF f)*(\FF g)$
for sufficiently regular tempered distributions 
$f$ and $g$, e.g., for $f, g \in \Ld$.
\vskip 0.2cm \noindent
\textbf{Proof of Prop.\ref{disug}.} We put for brevity 
$F := \FF f, ~G := \FF g$. By the relations between $\FF$, the $L^2$ norm 
and the convolution, we see that the thesis \rref{apart} is proved 
if we show that
\beq {\sqrt{1 + | \k |^2}~}^n | F \ast G | \in \Ld \qquad \mbox{and} 
\label{ispro} \feq
$$ \| {\sqrt{1 + | \k |^2}~}^n | F \ast G |~ \|_{L^2}  \leq 
\sum_{\ell \in \Lambda(n)}
\left( \barray{cc} n \\ \ell \farray 
\right)_{+} \| \left( {\sqrt{1 + | \k |^2}}^{~\ell}  | F | \right) \ast 
\left( {\sqrt{1 + | \k |^2}}^{~n - \ell}  | G | \right) \|_{L^2} $$
(note that the convolutions in the r.h.s are $L^2$ by our assumptions).
In order to derive Eq.\rref{ispro} we observe that, for $k \in \reali^d$, 
$$ {\sqrt{1 + | k |^2}~}^n | F \ast G | (k) \leq 
\int_{\reali^d} d h~{\sqrt{1 + | k |^2}~}^n~| F(k - h) |~ | G(h) | \leq
$$
$$ \leq \int_{\reali^d} d h~
\left(\sqrt{1 + | k - h |^2} + \sqrt{1 + | h |^2} \right)^n 
| F(k - h) |~ | G(h) | \leq $$
$$ \leq \sum_{\ell \in \Lambda(n)}
\left( \barray{cc} n \\ \ell \farray \right)_{+}
\int_{\reali^d} d h~{\sqrt{1 + | k - h |^2}~}^{\ell}~| F(k - h) |~ 
{\sqrt{1 + | h |^2}~}^{n - \ell}~| G(h) |~; $$
in the last two steps, we have employed Lemmas \ref{vettori} and 
\ref{newton}. Summing up, we have the pointwise inequality
\beq {\sqrt{1 + | \k |^2}~}^n | F \ast G |  \leq 
\sum_{\ell \in \Lambda(n)}
\left( \barray{cc} n \\ \ell \farray 
\right)_{+} \left( {\sqrt{1 + | \k |^2}}^{~\ell}  | F | \right) \ast
\left( {\sqrt{1 + | \k |^2}}^{~n - \ell}  | G | \right)~. \feq
The functions in the r.h.s. are $L^2$, so the same happens for
the l.h.s.;  taking the $L^2$ norms of both sides, we get the desired 
inequality in \rref{ispro}. 
\fine
\section{Inequalities for products
$\mbox{\boldmath $\DD_{\ell}(f) \DD_{m}(g)$}$.}
\label{deriva}
Our approach is similar to the usual argument for proving 
the Moser calculus inequality \cite{Mos} \cite{Maj} \cite{Zei},
with the following differences: 
the $\DD_{\ell}$ operators replace systematically 
the partial derivatives appearing in the cited works, and 
all the constants are estimated. In the sequel, 
$d \in \naturali_0$ is an arbitrary space dimension. \parn
\begin{lemma}
\label{prod1} 
\textbf{Lemma.} Let $\a, \ell, m$ be real numbers such that 
$\a > {d/2}$, $0 \leq \ell \leq \a$,~$m \geq 0$. 
For each $f \in \Ha$, 
$g \in H^{\ell + m}(\reali^d, \complessi)$ it is
$\DD_{\ell}(f) \DD_{m}(g) \in \Ld$, and
\beq \| \DD_{\ell} (f) \DD_{m}(g) \|_{L^2} \leq
E_{\ell, \a, d} S_{\a, d}~\| f \|_{\a}~ \| g \|_{\ell + m}~, \label{eprod1} \feq
with $S_{\a, d}$, 
$E_{\ell, a, d}$ defined as in Eq.s \rref{sad} \rref{elad}.
\end{lemma}
\textbf{Proof.} We put
$$ p := {2 \a \over \ell}~, 
\qquad q := {2 \a \over \a - \ell}~ $$
(intending $1/0 := \infty$); 
by construction $p, q \geq 2$ and $1/p + 1/q = 1/2$. 
By the H\"older and the imbedding inequality 
(Prop.\ref{sob}), the functions we consider are in the 
spaces indicated below, and 
$$ \| \DD_{\ell}(f) \DD_m(g) \|_{L^2} \leq 
\| \DD_{\ell}(f) \|_{L^p} \| \DD_m(g) \|_{L^q} \leq $$ 
$$ \leq
\left( S_{p, \a - \ell, d}~ 
\| \DD_{\ell}(f) \|_{\a - \ell } \right)~
\left( S_{q, \ell, d}~ \| \DD_m(g) \|_{\ell} \right) 
= S_{p, \a - \ell, d}~ S_{q, \ell, d}~
\| f \|_{\a}~ \| g \|_{\ell + m}~. \label{wit} $$
On the other hand, using Eq.s \rref{srd}\rref{srdi} 
with the above definitions of $p$ and $q$, 
and comparing with Eq.s
\rref{sad} \rref{elad}, one checks the equality
$S_{p, \a - \ell, d}~ S_{q, \ell, d} = E_{\ell, \a, d} S_{\a, d}$, 
yielding the thesis. \fine
\begin{lemma}
\label{prod2}
\textbf{Lemma.} Let $\a, \ell, m$ be real numbers such that 
$a > {d/2}$, $\ell, m \geq {a/2}$. For each
$f, g \in H^{\ell + m}(\reali^d, \complessi)$ 
it is $\DD_{\ell}(f) \DD_{m}(g) \in  \Ld$ and
\beq \| \DD_{\ell}(f) \DD_m(g) \|_{L^2} \leq
\left(16/27\right)^{d/4}~S_{\a, d}~\mbox{Max} 
\Big(\| f \|_{\a}~ \| g \|_{\ell + m}~,
\| f \|_{\ell + m}~\| g \|_{\a}\Big)~, 
\label{eprod2} 
\feq
where $S_{\a, d}$ is again as in Eq.\rref{sad}. 
\end{lemma}
\textbf{Proof.} We use the relation $1/4 + 1/4 = 1/2$ with the H\"older,  
the imbedding and the multiplicative interpolation inequality 
(Prop. \ref{segmento}). These ensure that the functions in consideration 
are in the spaces indicated below, and give the estimates
$$ \| \DD_{\ell}(f) \DD_{m}(g) \|_{L^2} \leq
\| \DD_{\ell}(f) \|_{L^4} \| \DD_m(g) \|_{L^4} 
\leq \left( S_{4,\a/2,d}\right)^2 
\| \DD_{\ell}(f) \|_{\a/2} \| \DD_{m}(g) \|_{\a/2} = $$
$$ = 
\left( S_{4,\a/2,d}\right)^2 \| f \|_{\ell+\a/2} ~\| g \|_{m+\a/2} \leq
\left( S_{4,\a/2,d}\right)^2 
\mbox{Max} \Big(\| f \|_{\a}~ \| g \|_{\ell + m}~, 
\| f \|_{\ell + m}~ \| g \|_{\a} \Big)~.
$$
On the other hand, comparing the definitions \rref{srd} 
and \rref{sad} \rref{elad} one checks that
$\left(S_{4, \a/2, d}\right)^2 = \left(16/27\right)^{d/4}~S_{\a, d}$
(this also equals $E_{\a/2, \a, d} S_{\a, d}$). The proof is 
concluded.
\fine
\section{Proof of Prop.\ref{upper}: upper bounds 
on $\mbox{\boldmath $K_{n, \a, d}$}$.}
\label{above}
Our approach will rely on 
the Leibnitz inequality (Prop.\ref{disug}) and
Lemmas \ref{prod1}, \ref{prod2}; by the way, the argument employed 
to derive the upper bounds of Prop.\ref{upper} on the constants 
$K_{n, \a, d}$ will also give a non conventional proof of Prop.\ref{nm}. 
We divide the proof in some steps.
\vskip 0.2cm \noindent
\textbf{Proof of Eq.s \rref{ub} \rref{knelad1}, 
case} $\mbox{\boldmath $~0 \leq n \leq d/2 < \a$}$~.
Prop.\ref{disug} reduces the problem to analysing the
products $\DD_{\ell}(f) \DD_{n - \ell}(g) \in \Ld$ 
for $\ell \in \Lambda(n)$ $\subset [0, n]$. To estimate them, 
we use Lemma \ref{prod1} (with $m = n - \ell$); this suffices to get the 
thesis. \fine
\textbf{Proof of Eq.s \rref{ub} \rref{knelad1}, 
case} $\mbox{\boldmath $~n \geq \a > d/2$}$~.
Again, we must analyse the
products $\DD_{\ell}(f) \DD_{n - \ell}(g)$ 
for $\ell \in \Lambda(n) \subset [0, n]$. \parn
If $\ell < \a/2$, we use Lemma \ref{prod1} (with $m = n - \ell$). 
If $\a/2 \leq \ell \leq n - \a/2$, we use Lemma \ref{prod2} 
(with $m = n - \ell$). Finally, if $n - \a/2 < \ell$, we employ Lemma 
\ref{prod1} with the pairs $(\ell, m)$, $(f, g)$ replaced by 
$(n - \ell, \ell)$, $(g, f)$. This yields the thesis.
\fine
\textbf{Proof of the weaker bounds \rref{wub1} \rref{wub2}.}
Eq.\rref{wub1} follows trivially from 
\rref{duen} \rref{ub} \rref{knelad1} and the fact that
$E_{n, \ell, a, d} \leq 1$ for all $\ell$. \parn
Let us prove Eq.\rref{wub2}, assuming $n \geq a > d/2$.
In this case, from $E_{\ell, \a, d} \leq 1$ we infer
$$ K_{n, \a, d} \leq~ 
S_{\a, d}~\left( 
\sum_{\ell < \a/2} 
\left( \barray{cc} n \\ \ell \farray 
\right)_{+} + \left(16/27\right)^{d/4}
\sum_{a/2 \leq \ell \leq n - \a/2} 
\left( \barray{cc} n \\ \ell \farray \right)_{+} + 
\sum_{\ell > n - \a/2}   
\left( \barray{cc} n \\ \ell \farray \right)_{+}~\right)~; $$
in all sums, it is intended that $\ell$ takes values in the lattice 
$\Lambda(n)$. The first and the last sum are 
equal to $\displaystyle{\sum_{j=0}^{a_n -1} \left( \barray{cc} n_{+} 
\\ j \farray \right)}$, with $a_n$ as in Eq. \rref{defna}; this implies
\beq K_{n, \a, d} \leq~ 
S_{\a, d}~\left(~ \left(16/27\right)^{d/4}~2^{n_{+}} 
+ 2 \sum_{j=0}^{a_n-1} \left( \barray{cc} n_{+} \\ j \farray \right)~
\left(1 - \left(16/27\right)^{d/4}~\right)~\right). 
\label{lab} \feq
To go, on we need the elementary inequality 
$\displaystyle{ \left( \barray{cc} m \\ j \farray \right) \leq 
\left( \barray{cc} m \\ m - k \farray \right) 
\left( \barray{cc} k \\ j \farray \right) }$, holding 
for $0 \leq j \leq k \leq m$~integers. Applying it with
 $m=n_+$, $k = a_n - 1$ and summing over $j$ we obtain from
\rref{lab} the thesis \rref{wub2}. \fine
\section{"Ground level" lower bounds on $\mbox{\boldmath $K_{n, \a, d}$}$.} 
\label{ground}
The aim of this Section is to prove Prop.\ref{gen}. The major 
step consists in proving the following 
\begin{lemma}
\label{major}
\textbf{Lemma.} 
Let $d \in \naturali_0$ 
and either $0 \leq n \leq d/2 < \a$ or $n \geq \a > d/2$; then, the 
sharp constant in Eq.\rref{in1} or \rref{in2} is such that
\beq K_{n, a, d} \geq {| f(0) | \over \| f \|_a} \label{tesi} \feq 
for each nonzero $f \in \Ha$. 
\end{lemma} 
\textbf{Remark.} Evaluation of $f$ at zero makes sense because
$\Ha \subset C(\reali^d, \complessi)$. 
\vskip 0.1cm \noindent
\textbf{Proof.} The idea is very simple; let us introduce it 
heuristically, say for $n \leq d/2 < a$. Let $f \in \Ha$, and put 
$g := \delta$ (the Dirac distribution). Formally, we have
$f g = f(0) \delta$,~ $\| f g \|_n = | f(0) |~ \| \delta \|_n$, 
$\| g \|_n = \| \delta \|_n $; 
inserting these trial functions into
Eq.\rref{ofc1}, and simplifying $\| \delta \|_n$ as if it were 
well defined, we get \rref{tesi}. All the rest of the proof is 
simply a rigorization of this idea, also working for $n \geq a > d/2$. \parn
First of all, we note that it suffices to prove the thesis \rref{tesi}
for 
\beq f = \FF^{-1} F~, \qquad F \in L^{\infty}(\reali^d, \complessi),~ 
F \neq 0,~\mbox{Supp $F$ ~bounded} \label{kind} \feq
(where Supp is the essential support; functions of the above kind 
are dense in all Sobolev spaces, and the evaluation map 
$f \mapsto f(0)$ is continuous in the $\|~\|_a$ norm). 
So, let us assume \rref{kind}, and consider a one parameter family of 
functions
\beq g_{\ep} := \FF^{-1} \left(G_{\ep}\right) \quad (\ep > 0)~; \qquad
G_{\ep}(k) := G(\ep k) \quad \mbox{for $k \in \reali^d$}~; \feq
$$ G \in C^1 (\reali^d, 
\complessi),~ G \neq 0,~ \mbox{Supp $G$~ bounded}~. $$
The function $g_{\ep}$ 
belongs to $H^n$ for any $n \geq 0$; we will ultimately consider 
the limit $\ep \vain 0$, under which $G_{\ep}$ behaves like a constant, and
$g_{\ep}$ like a constant $\times$ $\delta$. We have
$$ \| g_{\ep} \|^2_n = \int_{\reali^d} d k (1 + | k |^2)^n 
| G(\ep k) |^2 = {1 \over \ep^{2 n + d}} \int_{\reali^d} d k 
(\ep^2 + | k |^2)^n | G(k) |^2~, $$
whence
\beq \| g_{\ep} \|_n \barray{ccc} ~ \\ 
\sim \\ \scriptstyle{\ep \vain 0} \farray
{1 \over \ep^{n + d/2}}~\sqrt{\int_{\reali^d} | k |^{2 n} | G(k) |^2} = 
{1 \over \ep^{n + d/2}}~\|~ | \k |^n G \|_{L^2} \qquad \mbox{for  
each $n \geq 0$}~.  \label{as1} \feq
Due to this asymptotics, for $n > a$ the 
product $\| f \|_a \| g_{\ep} \|_n$ clearly dominates 
$\| f \|_n \| g_{\ep} \|_a$ when $\ep$ is small; so, we have 
\beq \mbox{Max}( \| f \|_a \| g_{\ep} \|_{n}, \| f \|_n \| g_{\ep} \|_{a} )
\barray{ccc} ~ \\ \sim \\ \scriptstyle{\ep \vain 0} \farray
{1 \over \ep^{n + d/2 }}~\| f \|_a~\| ~| \k |^n G \|_{L^2} 
\qquad \mbox{for $n \geq a$}~. \label{as2} \feq
Let us pass to evaluate $\| f g_{\ep} \|_n$, for any $n \geq 0$. 
Reexpressing the pointwise product via Fourier transform and 
convolution, we get
$$ \| f g_{\ep} \|^2_n = {1 \over (2 \pi)^d} \int_{\reali^d} d k 
(1 + | k |^2)^n  | (F \ast G_{\ep})(k) |^2 = $$
$$ = {1 \over (2 \pi)^d} 
\int_{\reali^d} d p \int_{\reali^d} d q ~\overline{F(p)} F(q) 
\int_{\reali^d} d k  (1 + | k |^2)^n~ 
\overline{G(\ep k - \ep p)} G(\ep k - \ep q) = $$
\beq = {1 \over (2 \pi)^d \ep^{2 n + d} } 
\int_{\reali^d} d p \int_{\reali^d} d q ~\overline{F(p)} F(q) 
\int_{\reali^d} d k  (\ep^2 + | k |^2)^n 
\overline{G(k - \ep p)} G(k - \ep q)~; \label{soth} \feq
in the last two equalities, we have explicitated the 
convolution and rescaled by $\ep$ the integration variable. 
With our assumptions on $G$, it is not difficult to prove that
$$ \mbox{last integral in Eq.\rref{soth}} 
\barray{ccc} ~ \\ \vain \\ \scriptstyle{\ep \vain 0} \farray
\int_{\reali^d} d p \int_{\reali^d} d q ~\overline{F(p)} F(q) 
\int_{\reali^d} d k  ~| k |^{2 n} | G(k) |^2 = $$
$$ = \left| \int_{\reali^d} d p ~F(p) \right|^2 \| | \k |^{n} G \|^2_{L^2}  = 
(2 \pi)^d  | f(0) |^2 ~\|~ | \k |^{n} G \|^{2}_{L^2}~. $$
Inserting this into Eq.\rref{soth}, we finally get
\beq {\| f g_{\ep} \|_n}  
\barray{ccc} ~ \\ \sim \\ \scriptstyle{\ep \vain 0} \farray
{1 \over \ep^{n + d/2} } | f(0) |~\| ~| \k |^n G \|_{L^2} \qquad \mbox{for 
all $n \geq 0$}~. \label{as3} \feq
Let $0 \leq n \leq d/2 < a$. We write Eq.\rref{ofc1} with $f$  
and $g = g_{\ep}$ as above; sending $\ep$ to zero, and using 
the asymptotics \rref{as3} \rref{as1}, we get the thesis \rref{tesi}. 
For $n \geq a > d/2$, the thesis \rref{tesi} 
follows from Eq.\rref{ofc2} and from the asymptotics \rref{as3} \rref{as2}. 
\fine
Now, we are ready to give the \vskip 0.1cm \noindent
\textbf{Proof of Prop.\ref{gen}:}
$\mbox{\boldmath $K_{n, \a, d} \geq S_{a, d}.$}$
Due to the previous Lemma, it suffices to find a nonzero function 
$f \in \Ha$ such that
\beq {| f(0) | \over \| f \|_a} = S_{a, d}~. \feq
As observed in \cite{art}, this equality holds for 
$f := f_{a, d} =$ 
the function in Eq.s \rref{fnd} \rref{gemac}, with $n=a$. \fine
\vskip 0.1cm \noindent
\textbf{Applications.} As previously noted, Prop.\ref{gen} combined 
with the upper bound \rref{ub} implies Corollary \ref{cogen}, i.e., 
the equality $K_{0, a, d} = S_{a, d}$.  Also, the combination of these results 
confines $K_{n, a, d}$ to a fairly small interval when $n$ is low 
(see, e.g., Eq.\rref{cogen1}). For $d=1,2,3$, 
$a = [d/2] + 1$ and $n < a$ integer, we find
\beq K_{0,1,1} = {1 \over \sqrt{2}}~, \feq 
\beq K_{0,2,2} = {1 \over 2 \sqrt{\pi}} ~, \qquad 
{1 \over 2 \sqrt{\pi}} \leq K_{1,2,2} \leq {1 \over 2 \sqrt{\pi}} 
\left(1 + {4 \sqrt{3} \over 9}\right)~, \feq 
\beq K_{0,2,3} = {1 \over 2 \sqrt{2 \pi}}~,
\qquad {1 \over 2 \sqrt{2 \pi}} \leq K_{1,2,3} \leq 
{1 \over 2 \sqrt{2 \pi}} \left(1 + {8 \over 3^{9/4}} \right)~. \feq
These results correspond to the numerical bounds already written in 
Sect.\ref{main}. For completeness, let us also give the functions 
$f = f_{a, d}$ for the above values of $a$ and $d$.
From the general representation \rref{gemac} in terms of the Macdonald functions
$K_{(~)}$ (and from the equality $\rho^{1/2}~ K_{1/2}(\rho) = 
\sqrt{\pi/2}~e^{-\rho})$ we get
\beq f_{1,1} = \sqrt{{\pi \over 2}}~ e^{-| \x |}, \qquad
f_{2,2} = {1 \over 2}~| \x |~ K_{1}(|\x|), \qquad 
f_{2,3} = \sqrt{{\pi \over 8}}~ e^{-| \x |}~. \label{forc} \feq
\section{"Bessel" lower bounds
on $\mbox{\boldmath $K_{n, \a, d}$}$.}
\label{bessel}
These bounds are expressed by Prop.\ref{lalem},
and rest on the functions $f_{\la,n, d}$ in Eq.\rref{lagemac}. 
In this section we will compute the norms of $f_{\la, n, d}$ and of its 
square $f_{\la,n, d}^2$; after this, in a number of cases we will explicitate
the bound they give on $K_{n, a, d}$
maximising w.r.t. $\la$. 
Our results will be expressed in terms of the Beta 
function $B(z, w) = \Gamma(z) \Gamma(w)/\Gamma(z + w)$ 
and of the Gauss hypergeometric function $F= {~}_{2} F_{1}$. 
\begin{lemma}
\label{norman}
\textbf{Lemma.} For $d \in \naturali_{0}$, $n > d/2$ and $\la > 0$ it is
$f_{\la, n, d} \in \Hn$, and
\beq \| f_{\la,n, d} \|^2_{n} = {2~ \pi^{d/2} \over \Gamma(d/2) \la^d}~
\int_{0}^{+\infty} d s~ s^{d-1} {(1 + \la^2 s^2)^n \over (1 + s^2)^{2 n} } = 
\label{firsteq} \feq
$$ = {\pi^{d/2} \over \Gamma(d/2) \la^d} 
\Big[ B(2 n - {d/2}, {d/2})~ F({d/2}, -n, 1 + {d/2} - 2 n; \la^2) + $$
$$ + \la^{4 n - d}~B(n - {d/2}, {d/2} - 2 n)~ 
F(2 n, n - {d/2}, 1 - {d/2} + 2 n; \la^2)~\Big]~. $$
In particular, for $n$ integer it is
\beq \| f_{\la,n, d} \|^2_{n} = {\pi^{d/2} \over \Gamma(d/2) \la^d}~
~\sum_{\ell=0}^n~
\left( \barray{c} n \\ \ell \farray \right)
B(\ell + {d/2}, 2 n - {d/2} - \ell)~\la^{2 \ell}~. 
\label{remar} \feq
\end{lemma}
\textbf{Proof.} We have
\beq \| f_{\la,n, d} \|^2_n = \int_{\reali^d} d k~(1 + | k |^2)^n 
| \FF f_{\la,n, d} |^2 = {1 \over \la^{2 d}}~ 
\int_{\reali^d} d k~{(1 + | k |^2)^n \over (1 + | k |^2/\la^2)^{2 n}}~. 
\label{thelast} \feq
On the other hand, on radially symmetric functions depending only on 
$\rho := | k |$, it is $\int_{\reali^d} d k~= 2 \pi^{d/2}/\Gamma(d/2)~
\int_{0}^{+\infty} d \rho~\rho^{d-1}$; from here, 
expressing the second integral of Eq.\rref{thelast} 
in terms of the rescaled radial 
variable $s = | k |/\la$ we get the first equality \rref{firsteq}. The 
second equality \rref{firsteq} comes from a known expression of the 
above integral over $s$ in terms of hypergeometric functions 
(in the singular cases $2 n - d/2 - 1 \in \naturali$, 
the first hypergeometric in \rref{firsteq} 
must be appropriately intended, as a limit from nonsingular values).
For $n$ integer, the integral over $s$ 
can be computed expanding the binomial 
$(1 + \la^2 s^2)^n$, and integrating term by term; when this is 
done, Eq. \rref{remar} follows recalling that 
$\int_{0}^{+\infty} 
d s~s^\alpha/(1 + s^2)^\gamma = (1/2) B(\alpha/2 + 1/2,
\gamma - \alpha/2 - 1/2)$. 
\fine 
\begin{lemma}
\label{normaa}
\textbf{Lemma.} For $d \in \naturali_{0}$, $n \geq a > d/2$ and $\la > 0$ it is
$f_{\la, n, d} \in \Ha$, and
\beq \| f_{\la, n, d} \|^2_{a} = {2~ \pi^{d/2} \over \Gamma(d/2) \la^d}~
\int_{0}^{+\infty} d s~ s^{d-1} {(1 + \la^2 s^2)^a \over (1 + s^2)^{2 n} } = 
\label{firsteqa} \feq
$$ = {\pi^{d/2} \over \Gamma(d/2) \la^d} 
\Big[ B(2 n - {d/2}, {d/2})~ F({d/2}, -a, 1 + {d/2} - 2 n; \la^2) + $$
$$ + \la^{4 n - d}~B(2 n - a - {d/2}, {d/2} - 2 n)~ 
F(2 n, 2 n - a - {d/2}, 1 - {d/2} + 2 n; \la^2)~\Big]~. $$
In particular, for $a$ integer it is
\beq \| f_{\la,n, d} \|^2_{a} = {\pi^{d/2} \over \Gamma(d/2) \la^d}~
~\sum_{\ell=0}^a~
\left( \barray{c} a \\ \ell \farray \right)
B(\ell + {d/2}, 2 n - {d/2} - \ell)~\la^{2 \ell}~. 
\label{remara} \feq
\end{lemma}
\textbf{Proof.} Proceed as in the proof of Lemma \ref{norman}. \fine
\begin{lemma}
\label{normaq}
\textbf{Lemma.} For $d \in \naturali_{0}$, $n > d/2$ and $\la > 0$ it is
$f_{\la, n, d}^2 \in \Hn$, and
\beq \| f_{\la, n, d}^2 \|^2_{n} 
= {2~ \pi^{d/2} \over \Gamma(d/2) \la^d}~
{\Gamma^2(2 n - d/2) \over \Gamma^2(2 n)}~\times \label{ppoo} \feq
$$ \times 
\int_{0}^{+\infty} d s~ s^{d-1} (1 + 4 \la^2 s^2)^n
F(2 n - {d/2}, n, n + {1/2}; - s^2)^2~. $$
\fine
\end{lemma}
\textbf{Remark.} If one is able to express 
the hypergeometric in the above equation in terms of 
elementary functions, for some integer
$n=n_0$, one can derive an expression via elementary functions
for all integers $n \geq n_0$, by repeated application 
of differentiation operations (see, e.g., 
\cite{Abr}). For $d$ odd and $n$ integer it is
$$ F(2 n - {d/2}, n, n + {1/2}; - s^2) = 
{\mbox{a polynomial of order ($n - d/2 - 1/2$) 
in $s^2$} \over (1 + s^2)^{2 n - d/2 - 1/2}}~; $$
inserting this into Eq.\rref{ppoo}, 
we can reduce the integral therein to a 
linear combination of elementary integrals of the form 
$\int_{0}^{+\infty} 
d s~s^\alpha/(1 + s^2)^\gamma$. \fine
\vskip 0.2cm \noindent
\textbf{Proof of Lemma \ref{normaq}.} It is $f_{\la,n,d}^2 \in \Hn$ by 
Prop.\ref{nm} and the fact that $f_{\la,n,d} \in \Hn$.  
To compute the norm of this function, we start 
from its Fourier transform. For all (sufficiently regular) radially symmetric 
functions 
\beq f : \reali^d \vain \complessi~, \qquad 
f(x) = \varphi(r)~, \qquad r := | x |~ \feq
the Fourier transform $\FF f$ is also radially symmetric and 
given by \cite{Boc}
\beq (\FF f)(k) = {1 \over \rho^{d/2 - 1}}~
\int_{0}^{+\infty} d r~ r^{d/2} J_{d/2 - 1}(\rho r) \varphi(r)~,
\qquad \rho := | k |~, \feq
where $J_{(~)}$ are the Bessel functions of the first kind. \parn
In particular, for the radial function $f=f^2_{\la,n,d}$ we obtain
from Eq.s \rref{gemac} \rref{lagemac} 
$$ \left(\FF f^2_{\la,n,d}\right)(k) = 
{1 \over \rho^{d/2 - 1}}~
\int_{0}^{+\infty} d r~ r^{d/2} J_{d/2 - 1}(\rho r) 
{(\la r)^{2 n - d} \over 2^{2 n - 2} \Gamma^2(n)}~K^2_{n-d/2}(\la r) = $$
$$ = {2^{2 - 2 n} \over 
\Gamma^2(n) \la^{d/2 + 1} \rho^{d/2 - 1}}~
\int_{0}^{+\infty} d r~ r^{2 n - d/2} J_{d/2 - 1}
({\rho \over \la} r)~K^2_{n-d/2}(r) \label{obt} = $$
$$ = {2^{2 - 2 n} \over \Gamma^2(n) \la^{d/2 + 1} \rho^{d/2 - 1}} \times
{\sqrt{\pi} \over 2^{d/2}} \left({\rho \over \la}\right)^{d/2 - 1} 
{\Gamma(n) \Gamma(2 n - d/2) \over 2 \Gamma(n+1/2)}~ 
F(2 n - {d \over 2}, n, n + {1 \over 2}; - {1 \over 4} {\rho^2 \over \la^2}) = $$
\beq = {1\over 2^{d/2} \la^{d}}~
{\Gamma(2 n - d/2) \over \Gamma(2 n)}~
F(2 n - {d \over 2}, n, n + {1 \over 2}; - {1 \over 4} {\rho^2 \over \la^2})
\label{nonlab} \feq
with $\rho := | k |$. (In the last three steps: we have rescaled $r$ into
$r/\lambda$; we have used a known relation between integrals of Bessel
functions and $F$; we have employed the identity 
$\Gamma(n) \Gamma(n + 1/2) = 2^{1 - 2 n} \sqrt{\pi}~ \Gamma(2 n)$). 
Putting \rref{nonlab} into the Fourier representation \rref{repfur} 
of $\| f^{2}_{\la,n,d} \|_{n}$, and 
introducing the scaled radial variable $s := (1/2) | k |/\la = 
(1/2) (\rho/\la)$, we readily obtain the thesis \rref{ppoo}. \fine
At last, we have explicit expressions for the norms of $f_{\la,n,d}$ and 
its square, to be inserted into Eq.\rref{dare}; after doing this and 
maximising with respect to $\lambda$ we obtain the "Bessel" lower 
bound on $K_{n,a,d}$. We will compute this in three cases. 
\vskip 0.2cm \noindent
\textbf{Case } 
$\mbox{\boldmath $d=1$, $n=a=1$}$~. Eq.\rref{forc} (rescaled by 
$\lambda$) and Lemma \ref{normaa} give
\beq f_{\la,1,1} = \sqrt{\pi \over 2} e^{-\la | \x |}~,
\qquad \| f_{\la, 1,1} \|_{1} = \sqrt{ \pi \over 2}~ 
\sqrt{ \la + {1 \over \la}}~.
\feq
To compute the norm of $f^2_{\la, 1,1}$ we do not even need Lemma \ref{normaq} 
because $f^2_{\la, 1,1} = \sqrt{ \pi /2 }~ f_{2 \la, 1, 1}$,  which implies
\beq \| f^2_{\la, 1,1} \|_{1} = {\pi \over 2}~ \sqrt{ 2 \la + 
{1 \over 2 \la}}~. \feq
The Bessel lower bound for the present case is 
\beq K_{1,1,1} \geq 
\mbox{Sup}_{\lambda > 0}~ \KK_{1,1,1}(\la)~, \qquad 
\KK_{1,1,1}(\la) := 
{ \|  f^2_{\la, 1, 1} \|_1 \over \| f_{\la, 1, 1} \|_{1}^{~2} }~. \feq
The function $\KK_{1,1,1}(~)$ attains its absolute maximum at 
$\la = \sqrt{9 + \sqrt{97}}/(2 \sqrt{2}) \simeq 1.53$, which 
yields the lower bound reported in Eq.\rref{just}:
\beq K_{1,1,1} \geq \KK_{1,1,1} \left( 
{\sqrt{9 + \sqrt{97}} \over 2 \sqrt{2}} \right) > 0.84~. \feq
\vskip 0.1cm \noindent
\textbf{Case } 
$\mbox{\boldmath $d=2$, $n=a=2$}$~. Eq. \rref{forc} 
and Lemma \ref{normaa} give
\beq f_{\la,2,2} = {\la \over 2}~ | \x |~ K_{1}(\la | \x |)~, \qquad
\| f_{\la, 2, 2} \|_{2} = \sqrt{\pi \over 3} 
\sqrt{ \la^2 + 1 + {1 \over \la^2}}~. \feq
Concerning the square $f^2_{\la, 2, 2}$~, by Lemma \ref{normaq} we have
\beq \| f^2_{\la, 2, 2} \|_{2} = {\sqrt{2 \pi} \over 3 \la}~
\sqrt{\int_{0}^{+\infty} d s~
s (1 + 4 \la^2 s^2)^2 F(3,2,{5/2}; - s^2)^2}~; 
\label{numer} \feq
the corresponding hypergeometric is such that
\beq F(3,2,{5/2}; - s^2) = {3 (2 s^2 - 1) \over 16 s^2 (1 + s^2)^2} + 
{3 (1 + 4 s^2) \over 16 s^3 (1 + s^2)^{5/2}}~\mbox{ArcSinh(s)}~, \feq 
and the integral in Eq.\rref{numer} can be computed numerically. 
The Bessel lower bound is, in the present case,
\beq K_{2,2,2} \geq 
\mbox{Sup}_{\lambda > 0}~\KK_{2,2,2}(\la)~, \qquad 
\KK_{2,2,2}(\la) := 
{ \|  f^2_{\la, 2, 2} \|_{2} \over \| f_{\la, 2, 2} \|_{2}^{~2}}~. \feq
A numerical study of the function $\KK_{2,2,2}(~)$ 
shows that it attains its absolute maximum at a 
point close to $\la = 1.35$. In agreement with Eq.\rref{ju36}, we have
\beq K_{2,2,2} \geq \KK_{2,2,2}(1.35) > 0.36~. \feq
\vskip 0.1cm \noindent
\textbf{Case } 
$\mbox{\boldmath $d=3$, $n=a=2$}$~. Eq.\rref{forc} 
and Lemma \ref{normaa} give
\beq f_{\la,2,3} = \sqrt{\pi \over 8}~e^{-\la | \x |}~, \qquad
\| f_{\la, 2, 3} \|_{2} = {\pi \over \sqrt{8}}~ 
\sqrt{ 5 \la + {2 \over \la} + {1 \over \la^3}}~.
\feq 
Also, being $f^2_{\la, 2 ,3 } = \sqrt{\pi/8}~f_{2 \la, 2, 3}$~, we have
\beq \| f^2_{\la, 2, 3} \|_{2} = { \pi^{3/2} \over 8}~ 
\sqrt{ 10 \la + {1 \over \la} + {1 \over 8 \la^3}}~. \feq
The Bessel lower bound is 
\beq K_{2,2,3} \geq 
\mbox{Sup}_{\lambda > 0}~\KK_{2,2,3}(\la)~, \qquad 
\KK_{2,2,3}(\la) := 
{ \|  f^2_{\la, 2, 3} \|_{2} \over \| f_{\la, 2, 3} \|_{2}^{~2}}~. \feq
The function $\KK_{1,1,1}(~)$ attains its absolute maximum at a 
point close to $\la = 1.31$, and in agreement with Eq.\rref{ju24} we get
\beq K_{2,2,3} \geq \KK_{2,2,3}(1.31) > 0.24~. \feq
\fine
\section{"Fourier" lower bounds
on $\mbox{\boldmath $K_{n, \a, d}$}$.}
\label{below}
As anticipated, these are based on the trial functions
$f_{p, \sigma, d}(x) := 
e^{i p x_1}~e^{- (\sigma/ 2) | x |^2}$~ (with $p, \sigma >0$; 
see Eq.\rref{fg}).
Let us apply Eq.s \rref{ofc1} \rref{ofc2} 
with $f=g = f_{p, \sigma, d}$, taking into account that 
$\left(f_{p, \sigma, d}\right)^2 = f_{2 p, 2 \sigma, d}$; 
in this way we obtain
\beq K_{n, \a, d} \geq { \| f_{2 p, 2 \sigma, d} \|_n \over 
\| f_{p, \sigma, d} \|_n \| f_{p, \sigma, d} \|_{\a} } 
\qquad \qquad\mbox{for each $p, \sigma$}~, \label{obtain} \feq
both for $0 \leq n \leq d/2 < \a$ and for $n \geq a > d/2$.
We wish to infer from here the lower bounds of Prop.\ref{lower}; 
this result will follow from a number of Lemmas. First of all,
we will evaluate the norms of the 
trial functions and give upper and lower bounds for them, 
both interesting for $p$ large and $\sigma/p^2$ small. Then,
will insert these bounds in \rref{obtain} and get lower bounds 
for $K_{n,a,d}$, depending on $p, \sigma$. Finally, we will choose 
$p$ and $\sigma$ suitably, and obtain the lower bounds on $K_{n,a,d}$ 
of Prop. \ref{lower}.
\begin{lemma}
\label{gauss}
\textbf{Lemma.} For each $p, \si > 0$ and $n \geq 0$, 
the Fourier transform and the $n$-th norm of $f_{p, \sigma, d}$ are 
given by
\beq \left(\FF f_{p, \sigma, d}\right)(k) = 
{1 \over \sigma^{d/2}}~e^{- \left( (k_1 - p)^2 + {k_2}^2 + 
... + {k_d}^2 \right)/ (2 \sigma) } ~,  \label{fouf} \feq
\beq \| f_{p, \sigma, d} \|^{2}_{n} = {1 \over \sigma^d}~
\int_{\reali^d} d k~\left( 1 + (k_1 + p)^2 + {k_2}^2 + ... +{k_d}^2 
\right)^n e^{- | k |^2 / \sigma}~. \label{normaf} \feq
\end{lemma}
\textbf{Proof.} An elementary computation relying on 
$\displaystyle{\int_{-\infty}^{+\infty} d \xi~ e^{h \xi} 
e^{- \alpha \xi^2}  = \sqrt{ {\pi \over \alpha} }~
e^{h^2/(4 \alpha) }}$~~(for $h \in \complessi, \alpha > 0$).
\fine
\parn
\begin{lemma}
\label{lbnf}
\textbf{Lemma.} For each $p, \si > 0$ and $n \geq 1/2$, it is
\beq \| f_{p, \si, d} \|^{2}_n \geq \pi^{d/2}~{p^{2 n} \over \sigma^{d/2}}~. 
\label{lbd} \feq
\end{lemma}
\textbf{Proof.} Eq. \rref{normaf} implies
$$ \| f_{p, \si, d} \|^{2}_n \geq {1 \over \sigma^d}~\int_{\reali}
d k_1 |k_1 + p|^{2 n}~ e^{- {k_1}^2 / \sigma} 
\int_{\reali^{d-1}} 
d k_2 ... d k_d~e^{- {k_2}^2 / \sigma} ... ~e^{- {k_d}^2 / \sigma} = $$
\beq =  {\pi^{d/2 - 1/2} \over \sigma^{d/2 + 1/2}}~\int_{\reali}
d k_1 |k_1 + p|^{2 n}~ e^{- {k_1}^2 / \sigma} \label{imp} \feq 
(of course, the intermediate integral over $\reali^{d-1}$ is intended 
to be $1$ if $d=1$). 
On the other hand, 
$$ \int_{\reali}
d k_1 |k_1 + p|^{2 n}~ e^{- {k_1}^2 / \sigma} \geq 
\int_{-p}^{+\infty} d k_1 (k_1 + p)^{2 n}~ e^{- {k_1}^2 / \sigma} = $$
$$ = p^{2 n} \sqrt{\sigma} 
\int_{-p/\sqrt{\sigma}}^{+\infty} 
d t~ (1 + {\sqrt{\sigma} \over p} t)^{2 n}~ e^{- t^2} \geq 
p^{2 n} \sqrt{\sigma} \left( \int_{-p/\sqrt{\sigma}}^{+\infty} 
d t~ e^{-t^2} + 2 n 
{\sqrt{\sigma} \over p}~\int_{-p /\sqrt{\sigma}}^{\infty} d t~
t e^{-t^2} \right) $$
(in the last two steps: the variable change $t = k_1/\sqrt{\si}$ has 
been performed, and the Bernoulli inequality 
$(1 + u)^m \geq 1 + m~u$ for $m \geq 1$, $u > -1$ has been employed 
with $u = (\sqrt{\sigma}/p)~t$ and $m = 2 n$). Computing the above two integrals,
we get
\beq \int_{\reali}
d k_1 |k_1 + p|^{2 n}~ e^{- {k_1}^2 / \sigma} \geq 
p^{2 n} \sqrt{\pi \si}~ U_{n}({\sqrt{\si} \over p})~, 
\quad U_{n}(\xi) := {1 + \mbox{Erf}(1/\xi) \over 2} + {n \xi \over \sqrt{\pi}}
~e^{-1/\xi^2}~. \label{opop} \feq
Here, $\mbox{Erf}$ denotes as usually the error function. An elementary 
analysis shows that, for each $n \geq 1/2$, 
the function $U_{n}(~)$ is monotonically increasing on 
the domain $(0, +\infty)$; on the other 
hand $U_n(\xi) \vain 1$ for $\xi \vain 0^{+}$. So, 
$U_{n}(\xi) > 1$ for all $\xi > 0$; inserting this 
into Eq. \rref{opop}, and the result into \rref{imp} we get the thesis. \fine
\begin{lemma}
\label{ubnf}
\textbf{Lemma.} Let $p, \si > 0$ and $n \geq 0$ be such that
$n {\sigma/p^2} < 1$. Then
\beq \| f_{p, \si, d} \|^{2}_n \leq \pi^{d/2}~
{e^{\dd{ {1 \over 1 - n \sigma/p^2}~{n^2 \sigma \over p^2} + 
{n \over p^2}}} \over 
\left(1 - n \sigma / p^2 \right)^{d/2}}~
{p^{2 n} \over \sigma^{d/2}}~. \label{ubd} \feq
\end{lemma}
\textbf{Proof.} The elementary inequality 
$1 + u \leq e^{u}$ for $u \in \reali$ implies
\beq (v + w)^n \leq v^n ~e^{n w/v} \qquad \mbox{for $v > 0$, $v + w > 0$, 
$n \geq 0$}~. \label{apply} \feq
Applying \rref{apply} with $v := p^2$, $w := 1 + 2 k_1 p + | k |^2$, 
and inserting the outcome into Eq.\rref{normaf} we obtain
$$ \| f_{p, \sigma, d} \|^{2}_{n} \leq {p^{2 n} \over \sigma^d}~
e^{n/p^2}~\int_{-\infty}^{+\infty} 
d k_1~e^{- (1/\sigma - n/p^2) k_1^2 + 2 n k_1/p }~
\Pi_{i=2}^d~\int_{-\infty}^{+\infty} 
d k_i~e^{- (1/\sigma - n/p^2) k_i^2 } $$       
(intending the last product to be $1$ if $d=1$). Our 
assumptions on $n, p, \sigma$ 
ensure all the above Gaussian integrals to converge, and their
computation yields the thesis \rref{ubd}. \fine
\begin{lemma}
\label{ready}
\textbf{Lemma.} Let either $1/2 \leq n \leq d/2 < \a$ or
$n \geq a > d/2$ and $p, \sigma >0$ be such that 
$\mbox{Max}(n, \a)~{\sigma/p^2} < 1$.  Then
\beq K_{n, \a, d} \geq {1 \over (2 \pi)^{d/4}}~{ 
\left(1 - n \sigma / p^2 \right)^{d/4}
\left(1 - \a \sigma / p^2 \right)^{d/4} \over 
e^{\dd{ {1 \over 1 - n \sigma/p^2}~{n^2 \sigma \over 2 p^2} + 
{1 \over 1 - \a \sigma/p^2}~{\a^2 \sigma \over 2 p^2} + 
{n + \a \over 2 p^2} }}  }~
{\sigma^{d/4} \over p^{\a}}~2^n~. \label{depen} \feq
\end{lemma}
\textbf{Proof.} We apply Eq.\rref{obtain}, using 
the following estimates: Eq.\rref{lbd} for $\| f_{2 p, 2 \sigma, d} \|_n$, 
Eq.\rref{ubd} for $\| f_{p, \sigma, d} \|_n$ and 
Eq.\rref{ubd} (with $n$ replaced by $\a$) 
for $\| f_{p, \sigma, d} \|_{\a}$. The term $2^n$ in Eq.\rref{depen} 
appears because $(2 p)^{n} = 2^n p^{n}$~. \fine
Now, we use the freedom we have for 
the choice of $p$, $\sigma$ in the previous Lemma; of course, 
we would like to maximise the r.h.s. of Eq.\rref{depen}, or 
at least to go close to the maximum. The choice we will present 
is the result of a careful inspection of Eq.\rref{depen}, and 
approximates well the maximum for large $n$; it leads directly to the
\vskip 0.2cm \noindent
\textbf{Proof of Prop.\ref{lower}.} We apply 
Lemma \ref{ready} with
\beq p := {\sqrt{n + \a} \over \sqrt{\lambda}}~, 
\qquad \sigma := {\mu / \lambda \over n + \a}~, \qquad 
\lambda > 0~, \quad 0 < \mu < \a~. \label{ins1} \feq
In this way, after a tedious computation we get
\beq K_{n, \a, d} \geq {1 \over (2 \pi)^{d/4}}~\phi_{\a/2 - d/4}(\lambda)~
\phi_{d/4}(\mu)~ v_{\mu, n, \a, d}~ {2^n \over (n + \a)^{\a/2 + d/4}}~, 
\label{ted} \feq 
$$ \phi_{\alpha}(\xi) := \xi^{\alpha}~ e^{-\xi/2} \qquad 
\mbox{for}~\alpha, \xi > 0~, $$ 
$$ v_{\mu, n, \a, d} :=  
\left(1 - {\mu \over n + \a} + {\mu^2 \a n \over (n + \a)^4} \right)^{d/4}~
e^{\dd{ {(2 a - \mu) \mu n + \mu \a^2 
\over 2 (n + \a)^2 - 2 \mu n }  - 
{\mu \a^2 \over 2 (n + \a)^2 - 2 \mu \a} } }~. $$
By construction, it is $\lim_{n \vain +\infty} v_{\mu, n, \a, d} = 1$
for each $\mu$, so this factor becomes irrelevant for large $n$. 
Now, we choose $\lambda, \mu$ so as to maximise the factors 
$\phi_{\a/2 - d/4}(\lambda) \phi_{d/4}(\mu)$. The maximising 
values are $\lambda = a - d/2$, $\mu = d/2$; inserting them
into Eq.\rref{ted}, we finally get the 
lower bound \rref{lb} for $K_{n, \a, d}$ 
(the factor $v_{n, \a, d}$ appearing 
in \rref{lb} is just the present coefficient 
$v_{\mu, n, \a, d}$ with $\mu = d/2$).  \parn
To conclude, we must derive the weaker bound \rref{wlb} for 
$n \geq \a > d/2$; this follows readily from the expression \rref{vnad} of 
$v_{n, \a, d}$ and from the inequalities
$$ 1 - {d \over 2 (n + a)} + 
{d^2 a n \over 4(n + \a)^4} \geq 1 - {d \over 2 \a}~, 
\quad {(4 \a - d) d n + 2 d \a^2 
\over 8 (n + \a)^2 - 4 d n }  - 
{d \a^2 \over 4 (n + \a)^2 - 2 \a d} \geq 0~. $$
\fine
\vskip 0.5cm \noindent
\textbf{Acknowledgments.} We are grateful to D.R. Adams and R.A. Adams for
some bibliographical indications in a preliminary phase of this 
work; in particular, we acknowledge
D.R. Adams for pointing out to us reference \cite{Run}. \parn
This paper has been partially supported by Istituto Nazionale di
Alta Matematica, G.N.F.M., and M.U.R.S.T~. 

\end{document}